\documentclass{gtart}

\def\ifplaintex{\expandafter\ifx\csname documentclass\endcsname\relax}

\ifplaintex 
\else
\fi

\expandafter\ifx\csname beginpicture\endcsname\relax
\expandafter\ifx\csname documentclass\endcsname\relax
\input pictex \else
\input prepictex \input pictex \input postpictex \fi\fi

\def\gt{{\mathsurround=0pt\it $\cal G\mskip-2mu$eometry \&\ 
$\cal T\!\!$opology}}        

\def\gtp{{\mathsurround=0pt\it $\cal G\mskip-2mu$eometry \&\ 
$\cal T\!\!$opology $\cal P\!$ublications}}  

\def\lognumber#1{\def\thelognumber{#1}}
\def\volumenumber#1{\def\thevolumenumber{#1}}
\def\papernumber#1{\def\thepapernumber{#1}}
\def\volumeyear#1{\def\thevolumeyear{#1}}

\def\pagenumbers#1#2{\def\startpage{#1}\def\finishpage{#2}}
\def\published#1{\def\publishdate{#1}}
\def\proposed#1{\def\theproposer{#1}}
\def\seconded#1{\def\theseconders{#1}}
\def\received#1{\def\receiveddate{#1}}
\def\revised#1{\def\reviseddate{#1}}
\def\accepted#1{\def\accepteddate{#1}}

\def\asciiaddress#1{\def\theasciiaddress{#1}}
\def\asciiemail#1{\def\theasciiemail{#1}}
\long\def\asciiabstract#1{\long\def\theasciiabstract{#1}}

\let\\\par\let\thelognumber\relax
\let\thevolumenumber\relax\let\thepapernumber\relax
\let\thevolumeyear\relax\let\thesamplenumber\relax\let\startpage\relax
\let\finishpage\relax\let\publishdate\relax\let\receiveddate\relax
\let\reviseddate\relax\let\accepteddate\relax\let\theasciititle\relax
\let\theasciiauthors\relax\let\theasciiaddress\relax
\let\theasciiabstract\relax
\let\theasciiemail\relax\let\theshortauthors\relax\let\theshorttitle\relax

\long\def\maketitlep{   

\count0=\startpage

\gt\hfill      
\beginpicture
\setcoordinatesystem units <0.33truein, 0.33truein> point at 2.2 0.9
\setplotsymbol ({$\cal G$})
\plotsymbolspacing=9truept
\circulararc 315 degrees from 0 1 center at 0 0
\setplotsymbol ({$\cal T$})
\circulararc 315 degrees from 1 -1 center at 1 0
\endpicture
\break
{\small\ifx\thesamplenumber\relax 
Volume \else Sample
\fi\thevolumenumber\ (\thevolumeyear)
\startpage--\finishpage\nl
Published: \publishdate}
\vglue 0.5truein plus 0.4fil minus 0.1truein

{\parskip=0pt\leftskip 0pt plus 1fil\def\\{\par\smallskip}{\ifplaintex\large
\else\Large\fi\bf\thetitle}\par\medskip}   

\vglue 0pt plus 0.1fil 

{\parskip=0pt\leftskip 0pt plus 1fil\def\\{\par}{\sc\theauthors}
\par\medskip}

\vglue 0pt plus 0.1fil 

{\small\parskip=0pt\let\newline\\
{\leftskip 0pt plus 1fil\def\\{\par}{\sl\theaddress}\par}
\expandafter\ifx\theemail\relax    
\relax\else\vglue 5pt plus 0.02fil minus 2pt\def\\{\stdspace{\rm 
and}\stdspace} 
\cl{Email:\stdspace\tt\theemail}\fi
\ifx\theurl\relax                  
\relax\else\vglue 5pt plus 0.02fil minus 2pt\def\\{\stdspace{\rm 
and}\stdspace}
\cl{URL:\stdspace\tt\theurl}\fi\par}

\vglue 7pt plus 0.3fil minus 3pt

{\bf Abstract}
\vglue 5pt plus 0.1fil minus 2pt

\theabstract

\vglue 7pt plus 0.3fil minus 3pt

{\bf AMS Classification numbers}\quad Primary:\quad \theprimaryclass

Secondary:\quad \thesecondaryclass

\vglue 5pt plus 0.3fil minus 2pt

{\bf Keywords}\quad \thekeywords

\vglue 10pt plus 0.5fil minus 5pt

{\small  Proposed: \theproposer\hfill Received: \receiveddate\nl
Seconded: \theseconders\hfill 
\ifx\reviseddate\relax                         
Accepted: \accepteddate                        
\else
Revised: \reviseddate                          
\fi}
\eject
}       

\let\maketitlepage\maketitlep
\let\maketitle\maketitlepage

\font\phead=cmsl9 scaled 950
\font=cmsl9 scaled 1050
\font\pnum=cmbx10 scaled 913
\font\lnum=cmbx10 
\font\pfoot=cmsl9 scaled 950
\font=cmsl9 scaled 1050
\ifplaintex
\headline{\vbox to 0pt{\vskip -4.5mm\line{\small\phead\ifnum
\count0=\startpage ISSN 1364-0380 (on line)
1465-3060 (printed) \hfill {\pnum\folio}\else\ifodd\count0\def\\{ }%
\ifx\theshorttitle\relax\thetitle\else\theshorttitle\fi\hfill{\pnum\folio}
\else\def\\{ and }{\pnum\folio}\hfill\ifx\theshortauthors\relax\theauthors
\else\theshortauthors\fi\fi\fi}\vss}}
\footline{\vbox to 0pt{\vglue 0mm\line{\small\pfoot\ifnum\count0=\startpage
\copyright\ \gtp\hfill\else
\gt, Volume \thevolumenumber\ (\thevolumeyear)\hfill\fi}\vss
}}
\else
\makeatletter
\def\@oddhead{{\small\lhead\ifnum\count0=\startpage ISSN 1364-0380 (on line)
1465-3060 (printed) \hfill {\lnum\number\count0}\else\ifodd\count0
\def\\{ }\ifx\theshorttitle\relax \thetitle \else\theshorttitle\fi\hfill
{\lnum\number\count0}\else\def\\{ and }{\lnum\number\count0}
\hfill\ifx\theshortauthors\relax 
\theauthors\else\theshortauthors\fi\fi\fi}}\def\@evenhead{\@oddhead}
\def\@oddfoot{\small\lfoot\ifnum\count0=\startpage\copyright\ \gtp\hfill\else
\gt, Volume \thevolumenumber\ (\thevolumeyear)\hfill\fi}
\def\@evenfoot{\@oddfoot}
\makeatother
\fi

\newwrite\gtoutfile
\long\gdef\makeheadfile{  
{\def\\{, }\def\s{ }
\immediate\openout\gtoutfile head.xxx
\immediate\write\gtoutfile{To: math@arxiv.org}
\immediate\write\gtoutfile{Subject: put or rep NNNNN:pppp}
\immediate\write\gtoutfile{--text follows this line--}
\immediate\write\gtoutfile{Proxy-for: \ifx\theasciiauthors\relax
\theauthors\else\theasciiauthors\fi\s<\ifx\theasciiemail\relax\theemail\else\theasciiemail\fi>}
\immediate\write\gtoutfile{\noexpand\\}
\immediate\write\gtoutfile{Authors: \ifx\theasciiauthors\relax
\theauthors\else\theasciiauthors\fi}
{\def\\{ }\immediate\write\gtoutfile{Title: \ifx\theasciititle\relax
\thetitle\else\theasciititle\fi}}
\immediate\write\gtoutfile{Subj-class: GT or SG or MG etc}
\immediate\write\gtoutfile{MSC-class: \theprimaryclass\ifx\thesecondaryclass\relax\else, \thesecondaryclass\fi}
\immediate\write\gtoutfile{Journal-ref: Geom. Topol. \thevolumenumber
(\thevolumeyear) \startpage-\finishpage}
\immediate\write\gtoutfile{Comments: Published by Geometry and Topology at}
\immediate\write\gtoutfile{\s\s http://www.maths.warwick.ac.uk/gt/GTVol\thevolumenumber/paper\thepapernumber.abs.html}
\immediate\write\gtoutfile{\noexpand\\}
\immediate\write\gtoutfile{}
\ifx\theasciiabstract\relax
\immediate\write\gtoutfile{\theabstract}\else
\immediate\write\gtoutfile{\theasciiabstract}\fi
\immediate\write\gtoutfile{}
\immediate\write\gtoutfile{\noexpand\\}
\immediate\write\gtoutfile{}
\immediate\closeout\gtoutfile}}  

\def\maketitlepage{\maketitlep\makeheadfile}
\let\maketitle\maketitlepage

\def\ifplaintex{\expandafter\ifx\csname documentclass\endcsname\relax}

\ifplaintex 
\else
\fi

\expandafter\ifx\csname beginpicture\endcsname\relax
\expandafter\ifx\csname documentclass\endcsname\relax
\input pictex \else
\input prepictex \input pictex \input postpictex \fi\fi

\def\gt{{\mathsurround=0pt\it $\cal G\mskip-2mu$eometry \&\ 
$\cal T\!\!$opology}}        

\def\gtp{{\mathsurround=0pt\it $\cal G\mskip-2mu$eometry \&\ 
$\cal T\!\!$opology $\cal P\!$ublications}}  

\def\lognumber#1{\def\thelognumber{#1}}
\def\volumenumber#1{\def\thevolumenumber{#1}}
\def\papernumber#1{\def\thepapernumber{#1}}
\def\volumeyear#1{\def\thevolumeyear{#1}}

\def\pagenumbers#1#2{\def\startpage{#1}\def\finishpage{#2}}
\def\published#1{\def\publishdate{#1}}
\def\proposed#1{\def\theproposer{#1}}
\def\seconded#1{\def\theseconders{#1}}
\def\received#1{\def\receiveddate{#1}}
\def\revised#1{\def\reviseddate{#1}}
\def\accepted#1{\def\accepteddate{#1}}

\def\asciiaddress#1{\def\theasciiaddress{#1}}
\def\asciiemail#1{\def\theasciiemail{#1}}
\long\def\asciiabstract#1{\long\def\theasciiabstract{#1}}

\let\\\par\let\thelognumber\relax
\let\thevolumenumber\relax\let\thepapernumber\relax
\let\thevolumeyear\relax\let\thesamplenumber\relax\let\startpage\relax
\let\finishpage\relax\let\publishdate\relax\let\receiveddate\relax
\let\reviseddate\relax\let\accepteddate\relax\let\theasciititle\relax
\let\theasciiauthors\relax\let\theasciiaddress\relax
\let\theasciiabstract\relax
\let\theasciiemail\relax\let\theshortauthors\relax\let\theshorttitle\relax

\long\def\maketitlep{   

\count0=\startpage

\gt\hfill      
\beginpicture
\setcoordinatesystem units <0.33truein, 0.33truein> point at 2.2 0.9
\setplotsymbol ({$\cal G$})
\plotsymbolspacing=9truept
\circulararc 315 degrees from 0 1 center at 0 0
\setplotsymbol ({$\cal T$})
\circulararc 315 degrees from 1 -1 center at 1 0
\endpicture
\break
{\small\ifx\thesamplenumber\relax 
Volume \else Sample
\fi\thevolumenumber\ (\thevolumeyear)
\startpage--\finishpage\nl
Published: \publishdate}
\vglue 0.5truein plus 0.4fil minus 0.1truein

{\parskip=0pt\leftskip 0pt plus 1fil\def\\{\par\smallskip}{\ifplaintex\large
\else\Large\fi\bf\thetitle}\par\medskip}   

\vglue 0pt plus 0.1fil 

{\parskip=0pt\leftskip 0pt plus 1fil\def\\{\par}{\sc\theauthors}
\par\medskip}

\vglue 0pt plus 0.1fil 

{\small\parskip=0pt\let\newline\\
{\leftskip 0pt plus 1fil\def\\{\par}{\sl\theaddress}\par}
\expandafter\ifx\theemail\relax    
\relax\else\vglue 5pt plus 0.02fil minus 2pt\def\\{\stdspace{\rm 
and}\stdspace} 
\cl{Email:\stdspace\tt\theemail}\fi
\ifx\theurl\relax                  
\relax\else\vglue 5pt plus 0.02fil minus 2pt\def\\{\stdspace{\rm 
and}\stdspace}
\cl{URL:\stdspace\tt\theurl}\fi\par}

\vglue 7pt plus 0.3fil minus 3pt

{\bf Abstract}
\vglue 5pt plus 0.1fil minus 2pt

\theabstract

\vglue 7pt plus 0.3fil minus 3pt

{\bf AMS Classification numbers}\quad Primary:\quad \theprimaryclass

Secondary:\quad \thesecondaryclass

\vglue 5pt plus 0.3fil minus 2pt

{\bf Keywords}\quad \thekeywords

\vglue 10pt plus 0.5fil minus 5pt

{\small  Proposed: \theproposer\hfill Received: \receiveddate\nl
Seconded: \theseconders\hfill 
\ifx\reviseddate\relax                         
Accepted: \accepteddate                        
\else
Revised: \reviseddate                          
\fi}
\eject
}       

\let\maketitlepage\maketitlep
\let\maketitle\maketitlepage

\font\phead=cmsl9 scaled 950
\font=cmsl9 scaled 1050
\font\pnum=cmbx10 scaled 913
\font\lnum=cmbx10 
\font\pfoot=cmsl9 scaled 950
\font=cmsl9 scaled 1050
\ifplaintex
\headline{\vbox to 0pt{\vskip -4.5mm\line{\small\phead\ifnum
\count0=\startpage ISSN 1364-0380 (on line)
1465-3060 (printed) \hfill {\pnum\folio}\else\ifodd\count0\def\\{ }%
\ifx\theshorttitle\relax\thetitle\else\theshorttitle\fi\hfill{\pnum\folio}
\else\def\\{ and }{\pnum\folio}\hfill\ifx\theshortauthors\relax\theauthors
\else\theshortauthors\fi\fi\fi}\vss}}
\footline{\vbox to 0pt{\vglue 0mm\line{\small\pfoot\ifnum\count0=\startpage
\copyright\ \gtp\hfill\else
\gt, Volume \thevolumenumber\ (\thevolumeyear)\hfill\fi}\vss
}}
\else
\makeatletter
\def\@oddhead{{\small\lhead\ifnum\count0=\startpage ISSN 1364-0380 (on line)
1465-3060 (printed) \hfill {\lnum\number\count0}\else\ifodd\count0
\def\\{ }\ifx\theshorttitle\relax \thetitle \else\theshorttitle\fi\hfill
{\lnum\number\count0}\else\def\\{ and }{\lnum\number\count0}
\hfill\ifx\theshortauthors\relax 
\theauthors\else\theshortauthors\fi\fi\fi}}\def\@evenhead{\@oddhead}
\def\@oddfoot{\small\lfoot\ifnum\count0=\startpage\copyright\ \gtp\hfill\else
\gt, Volume \thevolumenumber\ (\thevolumeyear)\hfill\fi}
\def\@evenfoot{\@oddfoot}
\makeatother
\fi

\newwrite\gtoutfile
\long\gdef\makeheadfile{  
{\def\\{, }\def\s{ }
\immediate\openout\gtoutfile head.xxx
\immediate\write\gtoutfile{To: math@arxiv.org}
\immediate\write\gtoutfile{Subject: put or rep NNNNN:pppp}
\immediate\write\gtoutfile{--text follows this line--}
\immediate\write\gtoutfile{Proxy-for: \ifx\theasciiauthors\relax
\theauthors\else\theasciiauthors\fi\s<\ifx\theasciiemail\relax\theemail\else\theasciiemail\fi>}
\immediate\write\gtoutfile{\noexpand\\}
\immediate\write\gtoutfile{Authors: \ifx\theasciiauthors\relax
\theauthors\else\theasciiauthors\fi}
{\def\\{ }\immediate\write\gtoutfile{Title: \ifx\theasciititle\relax
\thetitle\else\theasciititle\fi}}
\immediate\write\gtoutfile{Subj-class: GT or SG or MG etc}
\immediate\write\gtoutfile{MSC-class: \theprimaryclass\ifx\thesecondaryclass\relax\else, \thesecondaryclass\fi}
\immediate\write\gtoutfile{Journal-ref: Geom. Topol. \thevolumenumber
(\thevolumeyear) \startpage-\finishpage}
\immediate\write\gtoutfile{Comments: Published by Geometry and Topology at}
\immediate\write\gtoutfile{\s\s http://www.maths.warwick.ac.uk/gt/GTVol\thevolumenumber/paper\thepapernumber.abs.html}
\immediate\write\gtoutfile{\noexpand\\}
\immediate\write\gtoutfile{}
\ifx\theasciiabstract\relax
\immediate\write\gtoutfile{\theabstract}\else
\immediate\write\gtoutfile{\theasciiabstract}\fi
\immediate\write\gtoutfile{}
\immediate\write\gtoutfile{\noexpand\\}
\immediate\write\gtoutfile{}
\immediate\closeout\gtoutfile}}  

\def\maketitlepage{\maketitlep\makeheadfile}
\let\maketitle\maketitlepage

\lognumber{155}
\volumenumber{6}\papernumber{4}\volumeyear{2002}
\pagenumbers{69}{89}
\received{15 December 2000}
\accepted{28 February 2002}
\revised{28 February 2002}
\proposed{Joan Birman}
\seconded{Dieter Kotschick, Steven Ferry}
\published{1 March 2002}

\usepackage{amsmath,amssymb,graphicx}
\def\color[rgb]#1{\relax}

\newtheorem{thm}{Theorem}
\newtheorem{lemma}[thm]{Lemma}
\newtheorem{cor}[thm]{Corollary}
\newtheorem{prop}[thm]{Proposition}

\theoremstyle{remark}

\newtheorem*{definition*}{Definition}
\newtheorem*{remark*}{Remark}

\def\pmf{{\cal PMF}}
\def\mf{{\cal MF}}
\def\R{{\mathbb R}}
\def\r{{\mathbb R_+}}
\def\Z{{\mathbb Z}}
\def\H{{\mathbb H}}

\def\pa{pseudo-Anosov\ }

\def\C{{\cal C}}
\def\tqh{{\widetilde{QH}}}

\begin{document}

\title{Bounded cohomology of subgroups of\\mapping class groups}
  \authors{Mladen Bestvina\\Koji Fujiwara}
\address{Mathematics Department, University of Utah\\155 
South 1400 East, JWB 233\\Salt Lake City, UT 84112, USA}

\secondaddress{Mathematics Institute, Tohoku University\\Sendai,
980-8578, Japan}

\asciiaddress{Mathematics Department, University of Utah\\155 South
1400 East, JWB 233\\Salt Lake City, UT 84112, USA\\and\\Mathematics
Institute, Tohoku University\\Sendai, 980-8578, Japan}

\email{bestvina@math.utah.edu}
\secondemail{fujiwara@math.tohoku.ac.jp}
\asciiemail{bestvina@math.utah.edu, fujiwara@math.tohoku.ac.jp}

\begin{abstract} 
We show that every subgroup of the mapping class group $\text{\it
MCG}(S)$ of a compact surface $S$ is either virtually abelian or it
has infinite dimensional second bounded cohomology. As an application,
we give another proof of the Farb--Kaimanovich--Masur rigidity theorem
that states that $\text{\it MCG}(S)$ does not contain a higher rank
lattice as a subgroup.\end{abstract}

\asciiabstract{We show that every subgroup of the mapping class
group MCG(S) of a compact surface S is either virtually abelian or it
has infinite dimensional second bounded cohomology. As an application,
we give another proof of the Farb-Kaimanovich-Masur rigidity theorem
that states that MCG(S) does not contain a higher rank lattice as a
subgroup.}


\maketitlepage
\section{Introduction}

When $G$ is a discrete group, a {\it quasi-homomorphism} on $G$ is a
function $h\co G\to\R$ such that
$$\Delta(h):=\sup_{\gamma_1,\gamma_2\in
G}|h(\gamma_1\gamma_2)-h(\gamma_1)-h(\gamma_2)|<\infty.$$ The number
$\Delta(h)$ is the {\it defect} of $h$. Let ${\cal V}(G)$ be the
vector space of all quasi-homomorphisms $G\to\R$. By $BDD(G)$ and
respectively 
$HOM(G)=H^1(G;\R)$ denote the subspaces of ${\cal V}(G)$ consisting of
bounded functions and respectively homomorphisms. Note that
$BDD(G)\cap HOM(G)=0$. We will be concerned with the quotient spaces 
$$QH(G)={\cal V}(G)/BDD(G)\qquad\text{and}$$ 
$$\tqh(G)={\cal V}(G)/(BDD(G)+HOM(G))\cong QH(G)/H^1(G;\R).$$
There is an exact sequence
$$0\to H^1(G;\R)\to QH(G)\to H^2_b(G;\R)\to H^2(G;\R)$$ where
$H^2_b(G;\R)$ denotes the second bounded cohomology of $G$ (for the
background on bounded cohomology the reader is referred to
\cite{gromov:bc} and \cite{monod:book}). Since $\tqh(G)$ is the
quotient $QH(G)/H^1(G;\R)$ we see that $\tqh(G)$ can be also
identified with the kernel of $H^2_b(G;\R)\to H^2(G;\R)$. If $G\to G'$
is an epimorphism then the induced maps $QH(G')\to QH(G)$ and
$\tqh(G')\to\tqh(G)$ are injective.

Calculations of $\tqh(G)$ have been made for many groups $G$. In all such
cases $\tqh(G)$ is either 0 or infinite dimensional. $\tqh(G)$
vanishes when $G$ is amenable (see \cite{gromov:bc}) and also notably
when $G$ is a cocompact irreducible lattice in a semisimple Lie group
of real rank $>1$ \cite{bm:lattices}.

In the sequence of papers
\cite{epstein-fujiwara,fujiwara:plms,fujiwara:bs} the
second author has established a method for showing that $\tqh(G)$ is
infinite dimensional for groups $G$ acting on hyperbolic spaces and
satisfying certain additional conditions. This represents a
generalization of the argument of Brooks \cite{brooks:bc} that
$\dim\tqh(G)=\infty$ 
when $G$ is a nonabelian free group. Theorem \ref{blueprint} can
be viewed as a refinement of that method.

Not every group $G$ acting on a hyperbolic space has
$\dim\tqh(G)=\infty$. A nontrivial example is provided by an
irreducible cocompact lattice in $SL_2(\R)\times SL_2(\R)$ that acts
(discretely) on the product $\H^2\times\H^2$ of two hyperbolic
planes. Notice that the action given by projecting to a single factor
is highly non-discrete. Our contribution in this paper is to identify
what we believe to be the ``right'' condition on the action that
guarantees $\dim\tqh(G)=\infty$. The condition is termed $\text{\it WPD}$ (``weak
proper discontinuity'').

The main application is to the action of mapping class groups on curve
complexes. These were shown to be hyperbolic by Masur--Minsky
\cite{minsky-masur:cc1}. The action is far from discrete --- indeed,
the vertex stabilizers are infinite. However, we will show that $\text{\it WPD}$
holds for this action. As a consequence we will deduce the rigidity
theorem of Farb--Kaimanovich--Masur that mapping class groups don't
contain higher rank lattices as subgroups. More generally, if $\Gamma$
is not virtually abelian and $\dim\tqh(\Gamma)<\infty$ then $\Gamma$
does not occur as a subgroup of a mapping class group. In particular,
if the mapping class group $\text{\it MCG}(S)$ is not virtually abelian, then
$\dim\tqh(G)=\infty$. This settles Morita's Conjectures 6.19 and 6.21
\cite{morita:survey} 
in the affirmative.

We now proceed with a review of hyperbolic spaces and we introduce
some terminology needed in the paper.

When $X$ is a connected graph, we consider the path metric $d=d_X$ on $X$
by declaring that each edge has length 1. A {\it geodesic arc} is
a path whose length is equal to the distance between its endpoints. A
{\it bi-infinite geodesic} is a line in $X$ such that every finite
segment is geodesic. Recall \cite{mg:hyperbolic} that $X$ is said to
be {\it $\delta$--hyperbolic} if for any three geodesic arcs
$\alpha,\beta,\gamma$ in $X$ that form a triangle we have that
$\alpha$ is contained in the $\delta$--neighborhood of
$\beta\cup\gamma$. 

A map $\phi\co Y \to X$ from a metric space $Y$ is a {\it
$(K,L)$--quasi-isometric (qi) embedding} if $$\frac 1Kd_Y(y,y')-L\leq
d_X(\phi(y),\phi(y'))\leq Kd_Y(y,y')+L$$ for all $y,y'\in Y$. A {\it
$(K,L)$--quasi-geodesic} (or just quasi-geodesic when $(K,L)$ are
understood) is a $(K,L)$--qi embedding of an interval (finite or
infinite). A fundamental property of $\delta$--hyperbolic spaces is
that there is $B=B(K,L,\delta)$ such that any two finite
$(K,L)$--quasi-geodesics with common endpoints are within $B$ of each
other, and also any two bi-infinite quasi-geodesics that are finite
distance from each other are within $B$ of each other. A qi embedding
$Y\to X$ is a {\it quasi-isometry} if the distance between points of
$X$ and the image of the map is uniformly bounded.

An isometry $g$ of $X$ is {\it axial} if there is a bi-infinite
geodesic (called an {\it axis} of $g$) on which $g$ acts as a
nontrivial translation. Any axis of $g$ is contained in the
$2\delta$--neighborhood of any other axis of $g$. More generally, an
isometry $g$ of $X$ is {\it hyperbolic} if it admits an invariant
quasi-geodesic (we will refer to it as a {\it quasi-axis} or a
$(K,L)$--quasi-axis if we want to emphasize $K$ and $L$). We will often
blur the distinction between a quasi-axis and its image. There are
easy examples of hyperbolic isometries that are not axial, but whose
squares are axial (eg, take the ``infinite ladder'' consisting of two
parallel lines and rungs joining corresponding integer points, and the
isometry that interchanges the lines and moves rungs one unit
up). When the graph is allowed to be locally infinite, there are
similar examples of hyperbolic isometries none of whose powers are
axial. In our main application, the action of the mapping class group
on the curve complex, it is unknown whether powers of hyperbolic
elements are axial. We are thankful to Howie Masur and Yair Minsky for
bringing up this point. Note that any two $(K,L)$--quasi-axes of $g$ are
within $B(K,L,\delta)$ of each other.

Every quasi-axis of $g$ is oriented by the requirement that $g$ acts
as a positive translation. We call this orientation the {\it
$g$--orientation} of the quasi-axis. Of course, the
$g^{-1}$--orientation is the opposite of the $g$--orientation. More
generally, any sufficiently long $(K,L)$--quasi-geodesic arc $J$ inside
the $B(K,L,\delta)$--neighborhood of a $(K,L)$--quasi-axis
$\ell$ of $g$ has a natural orientation given by $g$: a point of
$\ell$ within $B(K,L,\delta)$ of the terminal endpoint of $J$
is ahead (with respect to the $g$--orientation of $\ell$) of a point of
$\ell$ within $B(K,L,\delta)$ of the initial endpoint of $J$.  
We call this
orientation of $J$ the {\it $g$--orientation}. We say that two quasi-geodesic
arcs are {\it $C$--close} if each is contained in the $C$--neighborhood
of the other, and we say that two oriented quasi-geodesic arcs are {\it
oriented $C$--close} if they are $C$--close and the distance between
their initial and also their terminal endpoints is $\leq C$.

\begin{definition*}
When $g_1$ and $g_2$ are hyperbolic elements of $G$ we will write
$$g_1\sim g_2$$ if for an arbitrarily long segment $J$ in a
$(K,L)$--quasi-axis for $g_1$ there is $g\in G$ such that $g(J)$ is
within $B(K,L,\delta)$ of a $(K,L)$--quasi-axis of $g_2$ and $g\co J\to
g(J)$ is orientation-preserving with respect to the $g_1$--orientation
on $J$ and the $g_2$--orientation on $g(J)$.
\end{definition*}
Replacing the constant $B(K,L,\delta)$ by a larger constant
would not change the concept since for every $C>0$ there is $C'>0$
such that for any $(K,L)$--quasi-geodesic arc $J$ contained in the
$C$--neighborhood of a $(K,L)$--quasi-geodesic $\ell$ it follows that
the arc obtained by removing the $C'$--neighborhood of each vertex is
contained in the $B(K,L,\delta)$--neighborhood of
$\ell$. Similarly,
the concept does not depend on the choice of $(K,L)$. In
particular:
\begin{itemize}
\item $\sim$ is an equivalence
relation.
\item  $g_1\sim g_2$ if and only if $g_1^k\sim g_2^l$ for any $k,l$ with
$kl>0$. 
\item If $g_1$ and $g_2$ have positive powers which are
conjugate in $G$ then $g_1\sim g_2$. 
\end{itemize}

Under our condition $\text{\it WPD}$ (see Section
\ref{section:wpd}) the converse of the third bullet also holds.

\begin{definition*}
We say that the action of $G$ on $X$ is {\it nonelementary} if there
exist at least two hyperbolic elements whose $(K,L)$--quasi-axes do not
contain rays within finite distance of each other (this distance can
be taken to be $B(K,L,\delta)$). The two hyperbolic elements are then
called {\it independent}.
\end{definition*}

\begin{thm}\label{blueprint}
Suppose a group $G$ acts on a $\delta$--hyperbolic graph $X$ by
isometries. Suppose also that the action is nonelementary and that
there exist independent hyperbolic elements $g_1,g_2\in G$ such that
$g_1\not\sim g_2$.

Then $\tqh(G)$ is infinite dimensional.
\end{thm}

\begin{remark*}
Special cases of this theorem are discussed in the earlier papers of
the second author:
\begin{itemize}
\item \cite{epstein-fujiwara}\qua $G$ is a word-hyperbolic group acting on
its Cayley graph,
\item \cite{fujiwara:plms}\qua the action of $G$ on $X$ is properly discontinuous,
\item \cite{fujiwara:bs}\qua $G$ is a graph of groups acting on the
associated Bass--Serre tree.
\end{itemize}
\end{remark*}

\begin{prop}\label{10}\qua
Under the hypotheses of Theorem \ref{blueprint} there is a sequence\break
$f_1,f_2,\cdots\in G$ of hyperbolic elements such that
\begin{itemize}
\item $f_i\not\sim f_i^{-1}$\qua for\qua $i=1,2,\cdots$,\qua and
\item $f_i\not\sim f_j^{\pm 1}$\qua for\qua $j<i$.
\end{itemize}
\end{prop}

\begin{proof}
Since $g_1$ and $g_2$ are independent, we may replace $g_1,g_2$ by high
positive powers of conjugates to ensure that the subgroup $F$ of $G$
generated by $g_1,g_2$ is free with basis $\{g_1,g_2\}$, each
nontrivial element of $F$ is hyperbolic, and $F$ is
quasi-convex with respect to the action on $X$ (see \cite[Section
5.3]{mg:hyperbolic}). We will call such free subgroups {\it Schottky
groups}. 
Let $T$ be the Cayley graph of $F$ with respect to the generating set
$\{g_1,g_2\}$. Then $T$ is a tree and each oriented edge has a label
$g_i^{\pm 1}$. Choose a basepoint $x_0\in X$ and construct an
$F$--equivariant map $\Phi\co T\to X$ that sends $1$ to $x_0$ and sends
each edge to a geodesic arc. Quasi-convexity implies that $\Phi$
is a $(K,L)$--quasi-isometric embedding for some $(K,L)$ and in particular
for every $1\neq f\in F$ the $\Phi$--image of the axis of $f$ in
$T$ is a $(K,L)$--quasi-axis of $f$. By $\ell_i$ denote the axis of
$g_i$ in $T$, $i=1,2$.

Choose positive constants $$0\ll n_1\ll m_1\ll k_1\ll l_1\ll n_2\ll
m_2\ll\cdots$$ and define $$f_i=g_1^{n_i}g_2^{m_i}g_1^{k_i}g_2^{-l_i}$$ for
$i=1,2,3,\cdots$. 

{\bf Claim 1}\qua $f_1\not\sim f_2$. 

The key to the proof is the
following observation. If $K',L',C$ are fixed and the exponents
$n_1,m_1,\cdots,k_2,l_2$ are chosen suitably large, then for any
sufficiently long $f_2$--oriented segment $S$ in the axis
$\ell_2\subset T$ of $f_2$ and any orientation preserving $(K',L')$--qi
embedding $\phi\co S\to\ell_1$ (with respect to the $f_1$--orientation of
$\ell_1$) there is a subsegment $S'$ of $S$ containing a string of
$\geq C$ edges labeled $g_2$ whose image (pulled tight) is a segment
in $\ell_1$ consisting only of edges labeled $g_1$ (this is true
because $m_2\gg n_1+m_1+k_1+l_1$ so the image will contain a whole
fundamental domain for the action of $f_1$ on $\ell_1$). Figure
\ref{fig1} illustrates the situation.

\begin{figure}[ht!]\label{fig1}
\begin{center}
\begin{picture}(0,0)%
\includegraphics{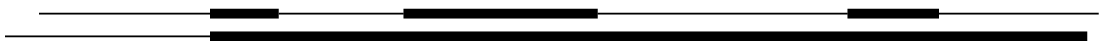}%
\end{picture}%
\setlength{\unitlength}{3947sp}%
\begingroup\makeatletter\ifx\SetFigFont\undefined%
\gdef\SetFigFont#1#2#3#4#5{%
  \reset@font\fontsize{#1}{#2pt}%
  \fontfamily{#3}\fontseries{#4}\fontshape{#5}%
  \selectfont}%
\fi\endgroup%
\begin{picture}(5272,604)(5014,-2441)
\put(7048,-1958){\makebox(0,0)[lb]{\smash{\SetFigFont{10}{12.0}{\rmdefault}{\mddefault}{\updefault}{\color[rgb]{0,0,0}$g_1^{k_1}$}%
}}}
\put(8195,-1958){\makebox(0,0)[lb]{\smash{\SetFigFont{10}{12.0}{\rmdefault}{\mddefault}{\updefault}{\color[rgb]{0,0,0}$g_2^{-l_1}$}%
}}}
\put(7704,-2396){\makebox(0,0)[lb]{\smash{\SetFigFont{10}{12.0}{\rmdefault}{\mddefault}{\updefault}{\color[rgb]{0,0,0}$g_1^{\pm n_2}$}%
}}}
\put(6010,-1958){\makebox(0,0)[lb]{\smash{\SetFigFont{10}{12.0}{\rmdefault}{\mddefault}{\updefault}{\color[rgb]{0,0,0}$g_1^{n_1}$}%
}}}
\put(6557,-1958){\makebox(0,0)[lb]{\smash{\SetFigFont{10}{12.0}{\rmdefault}{\mddefault}{\updefault}{\color[rgb]{0,0,0}$g_2^{m_1}$}%
}}}
\put(9202,-1958){\makebox(0,0)[lb]{\smash{\SetFigFont{10}{12.0}{\rmdefault}{\mddefault}{\updefault}{\color[rgb]{0,0,0}$g_1^{n_1}$}%
}}}
\end{picture}
\caption{Thick (thin)
lines represent strings of edges labeled $g_1$ ($g_2$).}
\end{center}
\end{figure}

Now assuming that $f_1\sim f_2$ let $I_2\subset\ell_2$ be a long arc,
let $J=\Phi(I_2)$, and let $g\in G$ be such that $g(J)$ is
$B(K,L,\delta)$--close to the $(K,L)$--quasi-axis $f_1(\ell_1)$ of $f_1$, with
matching orientations. Choose an arc $I_1\subset \ell_1$ so that
$\Phi(I_1)$ is $B(K,L,\delta)$--close to $g(J)$. Then there is
a $(K',L')$--quasi-isometry $I_2\to I_1$ obtained by composing 
$I_2\to \Phi(I_2)=J\to g(J)\to
\Phi(I_1)\to I_1$ and $(K',L')$ does not depend on the choices of
$n_1,m_1,\cdots,k_2,l_2$, only on $\delta$, $T$, and $\Phi$. Combining
this with the observation above, we conclude that $g$ takes a long
segment in an axis of a conjugate of $g_2$ uniformly close to an axis
of a conjugate of $g_1$ with matching orientation, contradicting the
assumption that $g_1\not\sim g_2$.

Similarly, $f_i\not\sim f_j$ for $i\neq j$.

{\bf Claim 2}\qua $f_1\not\sim f_2^{-1}$.

The proof is similar to the proof of Claim 1, only now one uses
$l_2\gg n_1+m_1+k_1+l_1$.

Similarly, $f_i\not\sim f_j^{-1}$ for $i\neq j$.

{\bf Claim 3}\qua If in addition $g_1\not\sim g_2^{-1}$ then $f_1\not\sim
f_1^{-1}$. 

If $f_1\sim
f_1^{-1}$, we obtain the situation pictured in Figure 2 where
a long string of $g_1$'s is close to a long string of $g_2$'s with
either the same or opposite orientation. Note that it is possible that
all such pairs of strings have opposite orientation so the assumption
$g_1\not\sim g_2^{-1}$ is necessary.
 
\begin{figure}[ht!]\label{fig2}
\begin{center}
\begin{picture}(0,0)%
\includegraphics{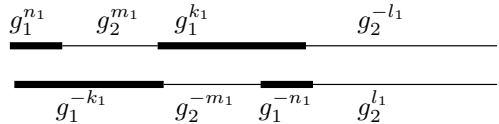}%
\end{picture}%
\setlength{\unitlength}{3947sp}%
\begingroup\makeatletter\ifx\SetFigFont\undefined%
\gdef\SetFigFont#1#2#3#4#5{%
  \reset@font\fontsize{#1}{#2pt}%
  \fontfamily{#3}\fontseries{#4}\fontshape{#5}%
  \selectfont}%
\fi\endgroup%
\begin{picture}(3116,684)(5966,-2510)
\put(7048,-1958){\makebox(0,0)[lb]{\smash{\SetFigFont{10}{12.0}{\rmdefault}{\mddefault}{\updefault}{\color[rgb]{0,0,0}$g_1^{k_1}$}%
}}}
\put(8195,-1958){\makebox(0,0)[lb]{\smash{\SetFigFont{10}{12.0}{\rmdefault}{\mddefault}{\updefault}{\color[rgb]{0,0,0}$g_2^{-l_1}$}%
}}}
\put(6010,-1958){\makebox(0,0)[lb]{\smash{\SetFigFont{10}{12.0}{\rmdefault}{\mddefault}{\updefault}{\color[rgb]{0,0,0}$g_1^{n_1}$}%
}}}
\put(6557,-1958){\makebox(0,0)[lb]{\smash{\SetFigFont{10}{12.0}{\rmdefault}{\mddefault}{\updefault}{\color[rgb]{0,0,0}$g_2^{m_1}$}%
}}}
\put(8195,-2500){\makebox(0,0)[lb]{\smash{\SetFigFont{10}{12.0}{\rmdefault}{\mddefault}{\updefault}{\color[rgb]{0,0,0}$g_2^{l_1}$}%
}}}
\put(7051,-2500){\makebox(0,0)[lb]{\smash{\SetFigFont{10}{12.0}{\rmdefault}{\mddefault}{\updefault}{\color[rgb]{0,0,0}$g_2^{-m_1}$}%
}}}
\put(7576,-2500){\makebox(0,0)[lb]{\smash{\SetFigFont{10}{12.0}{\rmdefault}{\mddefault}{\updefault}{\color[rgb]{0,0,0}$g_1^{-n_1}$}%
}}}
\put(6301,-2500){\makebox(0,0)[lb]{\smash{\SetFigFont{10}{12.0}{\rmdefault}{\mddefault}{\updefault}{\color[rgb]{0,0,0}$g_1^{-k_1}$}%
}}}
\end{picture}
\caption{Whenever thick and thin
lines are close, they are anti-parallel.}
\end{center}
\end{figure}
Similarly, if $g_1\not\sim g_2^{-1}$ then $f_i\not\sim f_i^{-1}$ for
all $i$.

We now finish the proof. If $g_1\not\sim g_2^{-1}$ then the above
claims conclude the argument. Otherwise, note that by Claims 1 and 2 we have
$f_1\not\sim f_2^{\pm 1}$. Now replace $(g_1,g_2)$ by $(f_1,f_2)$ and
repeat the construction.
\end{proof}

\section{Proof of Theorem \ref{blueprint}}

We will only give a sketch of the proof since it is a minor
generalization of results in
\cite{epstein-fujiwara,fujiwara:plms,fujiwara:bs} and
the proof uses the same techniques.

We start by recalling the basic construction of quasi-homomorphisms in
this setting. The model case of the free group is due to Brooks
\cite{brooks:bc}. 

Let $w$ be a finite (oriented) path in $X$. By $|w|$ denote the length
of $w$. For $g\in G$ the
composition $g\circ w$ is
a {\it copy} of $w$.  Obviously $|g
\circ w | = |w|$.

Let $\alpha$ be a finite path. We define 
$$|\alpha|_w = \{ \mbox{the maximal number of non-overlapping copies of
$w$ in $\alpha$} \}.$$

Suppose that $x,y \in X$ are two vertices and that 
$W$ is an integer with $0 < W < |w| $.
We define the integer
$$c_{w,W}(x,y) = d(x,y) - \inf_{\alpha} (|\alpha|-W|\alpha|_w),$$
where $\alpha$ ranges over all paths from $x$ to $y$.
Note that if $\alpha$ is such a path that contains a subpath whose
length is large compared to the distance between its endpoints, then
replacing this subpath by a geodesic arc between the endpoints
produces a new path with smaller $|\alpha|-W|\alpha|_w$. This
observation leads to the following lemma. 

\begin{lemma}{\rm \cite[Lemma 3.3]{fujiwara:plms}}\label{qg}\qua
Suppose a path $\beta$ realizes the infimum above. Then $\beta$ is a
$(\frac{|w|}{|w|-W},\frac{2W|w|}{|w|-W})$--quasi-geodesic.\qed
\end{lemma}

Replace $g_1,g_2$ by large positive powers if necessary, let $F$ be
the subgroup of $G$ generated by $g_1,g_2$, and let $\Phi\co T\to X$ be
an $F$--equivariant map with $\Phi(1)=x_0$ as in the proof of
Proposition \ref{10}. If $w\in F$ is cyclically reduced as a word in
$g_1,g_2$ (equivalently, if its axis passes through $1\in T$) then by
the quasi-convexity of $F$ in $G$ we have (see Figure 3)
$$d(x_0,w^n(x_0))\geq n(d(x_0,w(x_0))-2B)$$ where $B=B(K,L,\delta)>0$ is 
independent of $w$ and $n$. 

\begin{figure}[ht!]\label{fig3}
\begin{center}
\begin{picture}(0,0)%
\includegraphics{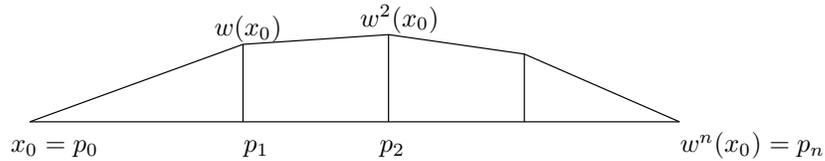}%
\end{picture}%
\setlength{\unitlength}{3947sp}%
\begingroup\makeatletter\ifx\SetFigFont\undefined%
\gdef\SetFigFont#1#2#3#4#5{%
  \reset@font\fontsize{#1}{#2pt}%
  \fontfamily{#3}\fontseries{#4}\fontshape{#5}%
  \selectfont}%
\fi\endgroup%
\begin{picture}(4216,998)(1126,-1922)
\put(2406,-1144){\makebox(0,0)[lb]{\smash{\SetFigFont{10}{12.0}{\rmdefault}{\mddefault}{\updefault}{\color[rgb]{0,0,0}$w(x_0)$}%
}}}
\put(2588,-1875){\makebox(0,0)[lb]{\smash{\SetFigFont{10}{12.0}{\rmdefault}{\mddefault}{\updefault}{\color[rgb]{0,0,0}$p_1$}%
}}}
\put(1126,-1875){\makebox(0,0)[lb]{\smash{\SetFigFont{10}{12.0}{\rmdefault}{\mddefault}{\updefault}{\color[rgb]{0,0,0}$x_0=p_0$}%
}}}
\put(3441,-1875){\makebox(0,0)[lb]{\smash{\SetFigFont{10}{12.0}{\rmdefault}{\mddefault}{\updefault}{\color[rgb]{0,0,0}$p_2$}%
}}}
\put(5330,-1875){\makebox(0,0)[lb]{\smash{\SetFigFont{10}{12.0}{\rmdefault}{\mddefault}{\updefault}{\color[rgb]{0,0,0}$w^n(x_0)=p_n$}%
}}}
\put(3320,-1083){\makebox(0,0)[lb]{\smash{\SetFigFont{10}{12.0}{\rmdefault}{\mddefault}{\updefault}{\color[rgb]{0,0,0}$w^2(x_0)$}%
}}}
\end{picture}
\caption{$nd(x_0,w(x_0))\leq
d(x_0,w^n(x_0))+d(p_1,w(x_0))+d(p_2,w^2(x_0))+\cdots 
\leq d(x_0,w^n(x_0))+2nB$}
\end{center}
\end{figure}

Following \cite{fujiwara:plms} we fix an integer
$W\geq 3B$ and will only consider $w$ with $|w|>W$. Thus, paths $\beta$
as in Lemma \ref{qg} will be quasi-geodesics with constants
independent of $w$ and the endpoints, and $\beta$ is contained in a
uniform neighborhood, say $D$--neighborhood, of any geodesic joining
the endpoints of $\beta$. We will also omit $W$ from notation and
write $c_w$. For every $f\in F$ choose a
geodesic $\gamma_f$ from $x_0$ to $f(x_0)$. We find it convenient to
denote the concatenation
$$\gamma_f f(\gamma_f) f^2(\gamma_f)\cdots f^{a-1}(\gamma_f)$$ by
$f^a$.

Define $h_w\co G\to\R$ by
$$h_{w}(g) =c_{w}(x_0,g(x_0)) - c_{w^{-1}}(x_0,g(x_0)).$$

\begin{prop}{\rm\cite[Proposition 3.10]{fujiwara:plms}}\label{qh}\qua
The map
$h_{w}\co G \to {\R}$ is a quasi-homomorphism. Moreover, the defect
$\Delta(h_w)$ is uniformly bounded independently of $w$.\qed
\end{prop}

\begin{prop}\label{13}
Suppose $1\neq f\in F$ is cyclically reduced and $f\not\sim f^{-1}$. Then
there is $a>0$ such that $h_{f^a}$ is unbounded on $<f>$. Moreover, if
$f^{\pm 1}\not\sim f'\in F$ then $h_{f^a}$ is 0 on $<f'>$ for sufficiently
large $a>0$.
\end{prop}

\begin{proof}
It is clear that $c_{f^a}$ is unbounded on $<f>$ for any $a>0$: If we
use $\alpha=(f^a)^n$ as a competitor path we have
$$c_{f^a}(f^{an})\geq d(x,f^{an}(x))-(n|f^a|-3Bn)\geq Bn.$$

If
$a>0$ is large, then there are no copies of $f^{-a}$ in the $D$--neighborhood
of an axis of $f$, which implies that $c_{f^{-a}}$ is zero on $<f>$.

The proof of the other claim is similar.
\end{proof}

\begin{proof}[Proof of Theorem \ref{blueprint}]
Let $f_1,f_2,\cdots$ be the sequence from Proposition \ref{10}.  We
assume in addition (without loss of generality) that each $f_i$ is
cyclically reduced. Define $h_i\co G\to\R$ as $h_i=h_{f^{a_i}}$ where
$a_i>0$ is chosen as in Proposition \ref{13} so that $h_i$ is
unbounded on $<f_i>$ and so that it is 0 on $<f_j>$ for $j<i$ (a high
power of $f_i$ cannot be translated into a $B$--neighborhood of an axis
of $f_j$). It follows that $[h_i]\in QH(G)$ is not a linear
combination of $[h_1],\cdots,[h_{i-1}]$, ie, the sequence $[h_i]$
consists of linearly independent elements. We can easily arrange that
$F$ is contained in the commutator subgroup of $G$ (this is automatic
if $g_1,g_2$ are in the commutator subgroup; otherwise, replace
$g_1,g_2$ by
$g_1'=g_1^Ng_2^Mg_1^{-M}g_2^{-N},g_2'=g_1^Kg_2^Lg_1^{-K}g_2^{-L}$ with
$0\ll N\ll M\ll K\ll L$ ---
as in Proposition \ref{10} it follows that $g_1'\not\sim g_2'$). In
that case any homomorphism $G\to\R$ vanishes on $F$ and it follows
that the sequence $[h_i]$ in $\tqh(G)$ consists of linearly
independent elements.
\end{proof}

\begin{remark*}
The argument shows that there is an embedding of $\ell^1$ into
$QH(G)$. If $t=(t_1,t_2,\cdots)\in \ell^1$ then $$h_t=\sum
t_ih_i\co G\to\R$$ is a well-defined function ($|f_i|\to\infty$ implies
that for any $g\in G$ only finitely many $h_i(g)$ are nonzero) and it
is a quasi-homomorphism (because the defect of $h_i$ is uniformly
bounded independently of $i$), and $t\neq 0$ implies $h_t$ is
unbounded (on $<f_i>$ where $i$ is the smallest index with $t_i\neq
0$). Similarly, one can argue that $\widetilde {QH}(G)$ contains
$\ell^1$.
\end{remark*}

\begin{remark*} Instead over $\R$ one can work over $\Z$ and consider
$H^2_b(G;\Z)$ and
quasi-homomorphisms $G\to\Z$. The
quasi-homomorphisms $h_w$ constructed above are integer-valued, and
therefore it follows that there are infinitely many linearly
independent elements in the kernel of $H^2_b(G;\R)\to H^2(G;\R)$ which
are in the image of $H^2_b(G;\Z)$.
\end{remark*}

\section{Weak Proper Discontinuity}\label{section:wpd}

\begin{definition*}
We say that the action of $G$ on $X$ satisfies $\text{\it WPD}$ if 
\begin{itemize}
\item $G$ is not
virtually cyclic,
\item $G$ contains at least one element that acts on $X$ as
a hyperbolic isometry, and 
\item for every hyperbolic element $g\in G$, every $x\in X$, and every
$C>0$ there exists $N>0$ such that
the set 
$$\{\gamma\in G|d(x,\gamma(x))\leq C,d(g^N(x),\gamma g^N(x))\leq C\}$$
is finite.
\end{itemize}
\end{definition*}

\begin{prop} \label{basic}
Suppose that $G$ and $X$ satisfy $\text{\it WPD}$. Then
\begin{itemize}
\item [\rm(1)]
for every hyperbolic $g\in G$ the centralizer $C(g)$ is virtually
cyclic,
\item [\rm(2)]
for every hyperbolic $g\in G$ and every $(K,L)$--quasi-axis $\ell$ for
$g$ there is a constant $M=M(g,K,L)$ such that if two translates $\ell_1$,
$\ell_2$ of $\ell$ contain (oriented) segments of length $>M$ that
are oriented $B(K,L,\delta)$--close then $\ell_1$ and $\ell_2$
are oriented $B(K,L,\delta)$--close and moreover the
corresponding conjugates $g_1,g_2$ of $g$ have positive powers which
are equal,
\item [\rm(3)] the action of $G$ on $X$ is nonelementary,
\item [\rm(4)] $f_1\sim f_2$ if and only if
some positive powers of $f_1$ and $f_2$ are conjugate,
\item [\rm(5)] there exist hyperbolic $g_1,g_2$ such that $g_1\not\sim g_2$.
\end{itemize}
\end{prop}

\begin{proof} (1)\qua We will show that $<g>$ has finite index in
$C(g)$. Let $f_1,f_2,\cdots$ be an infinite sequence of elements of
$C(g)$. Choose a $(K,L)$--quasi-axis $\ell$ for $g$ and let $x\in
\ell$. Since $g$ and $f_i$ commute, $f_i(\ell)$ is also a
$(K,L)$--quasi-axis for $g$ and the distance between $\ell$ and
$f_i(\ell)$ is uniformly bounded by $B(K,L,\delta)$. Let $k_i\in\Z$ be
such that $d(f_i(x),g^{k_i}(x))$ minimizes the distance between
$f_i(x)$ and the $g$--orbit of $x$. Thus $d(f_i(x),g^{k_i}(x))$ is
uniformly bounded (by $B(K,L,\delta)$ plus the diameter of the
fundamental domain for the action of $g$ on $\ell$); call such a bound
$C$. Let $N$ be from the definition of $\text{\it WPD}$. We note that
$f_ig^{-k_i}$ move both $x$ and $g^N(x)$ by $\leq C$. Therefore, the
set of such elements is finite. From $f_ig^{-k_i}=f_jg^{-k_j}$ we
conclude that $f_i$ and $f_j$ represent the same $<g>$--coset and the
claim is proved.

(2)\qua Denote by $g_1$ and $g_2$
the corresponding conjugates of $g$. For notational simplicity, we
will first assume that $g$ is axial and that $\ell$ is an axis of $g$.

Without loss of generality we assume $g_1=g$. Choose $x\in \ell$ and
let $N$ be as in the definition of $\text{\it WPD}$ for $g,x,C=4\delta$. Let
$P$ be the size of the finite set from the definition of $\text{\it WPD}$. If
$\ell_1=\ell$ and $\ell_2$ contain oriented $2\delta$--close arcs $J_1$
and $J_2$ of length $>(P+N+2)\tau_g$ ($\tau_g$ is the translation length
of $g$) then the elements $g_1^ig_2^{-i}$ move each point of the
terminal subarc of $J_2$ of length $(N+2)\tau_g$ a distance $\leq 4\delta$
for $i=0,1,\cdots,P$. It follows 
that $g_1^ig_2^{-i}=g_1^jg_2^{-j}$ for distinct $i,j$ so that
$g_1,g_2$ have equal positive powers.

In general, when $\ell$ is only a quasi-axis, one can
generalize the above paragraph by replacing $4\delta$ etc. by larger
constants that depend on $(K,L)$ and $\delta$. Alternatively, one can
modify $X$ to make $g$ axial: simply attach an infinite ladder (the
1--skeleton of an infinite strip) along one of the two infinite lines
to $\ell$; then attach such ladders equivariantly to obtain a
$G$--space. Finally, subdivide each rung and each edge in $X$
into a large number $Q$ of edges in order to arrange
that the ``free'' lines in the attached ladders are geodesics and axes
for the corresponding conjugates of $g$. The group $G$ continues to
act on the new space $X'$ which is quasi-isometric to $X$. The
statement for $X'$ implies the statement for $X$.

(3)\qua Let $g$ be a hyperbolic element. Again, without loss of generality,
 we will assume that $g$ has an axis $\ell$. We aim to show that some
 translate of $\ell$, say $h(\ell)$, has both ends distinct from
 $\ell$, since then $g$ and $hgh^{-1}$ are independent hyperbolic
 elements.

Suppose first that there is $h\in G$ such that $h(\ell)$ is not in the
 $2\delta$--neighborhood of $\ell$, but it is asymptotic in one
 direction, ie, a ray in $h(\ell)$ is contained in the
 $2\delta$--neighborhood of $\ell$. From (2) we see that $\ell$ and
 $h(\ell)$ cannot contain segments of length $>M(g)$ that are oriented
 $2\delta$--close to each other; in particular, one of $g$, $hgh^{-1}$
 moves towards the common end, say $\infty$, and the other moves away
 from it. Now consider $g^N(h(\ell))$ for large $N$. This is a
 bi-infinite geodesic with one end $\infty$ and the other end distinct
 from the ends of $\ell$, $h(\ell)$. The translates $h(\ell)$ and
 $g^N(h(\ell))$ violate (2) as they have oriented rays within $2\delta$
 of each other.

It remains to consider the case when every translate of $\ell$ is
within $2\delta$ of $\ell$. After passing to a subgroup of $G$ of
index 2 if necessary, we may assume that $G$ preserves the ends of
$\ell$. Now proceed as in (1) to show that $<g>$ has finite index in
$G$. 

(4)\qua This is similar to (2). We assume for simplicity that $f_1,f_2$ are
axial. By $\tau_k$ denote the translation length of $f_k$,
$k=1,2$. Let $N$ be as in the definition of $\text{\it WPD}$ for $g=f_1$ with
respect to some $x$ in an axis $\ell$ of $f_1$ and $C=4\delta$. Let
$P$ be the size of the corresponding finite set. Now assume that $f_2$
has an axis that admits a segment of length
$>(N+2)\tau_1+P\tau_1\tau_2$ which is oriented $2\delta$--close to a
segment in $\ell$. Consider $f_2^{i\tau_1}f_1^{-i\tau_2}$,
$i=0,1,\cdots,P$. As before, we conclude that
$f_2^{i\tau_1}f_1^{-i\tau_2}=f_2^{j\tau_1}f_1^{-j\tau_2}$ for some
$i\neq j$; thus $f_1$ and $f_2$ have common positive powers.

(5)\qua Since the action of $G$ on $X$ is nonelementary, we can choose a
   Schottky subgroup $F\subset G$. Let $1\neq f\in F$. For notational
   simplicity we will assume that all nontrivial elements of $F$ are
   axial and in fact that there is an $F$--invariant totally geodesic
   tree $T\subset X$ (this can be arranged by modifying $X$ as in the
   proof of (2) except that now one attaches the 1--skeleton of
   (tree)$\times I$ instead of (line)$\times I$). Let $\ell\subset T$
   be an axis of $f$. Then $\text{\it WPD}$ provides a segment $J$ in $\ell$ such
   that the set of $g\in G$ that move each point of $J$ by $\leq
   \tau_f+2\delta$ is finite.

Now consider an infinite sequence $f_1,f_2,\cdots$ of elements of $F$
with distinct (and hence non-parallel) axes
$\ell_1,\ell_2,\cdots\subset T$
that overlap $\ell$ in (oriented) finite intervals that contain
$J$. If $f_i\sim f$ then, according to (4), there is $g_i\in G$ such
that $g_i(\ell_i)$ is $2\delta$--close to $\ell$. Replacing $g_i$ by
$f^{a_i}g_i$ if necessary, we may assume that $g_i$ moves each point
of $J$ by $\leq \tau_f+2\delta$. Thus $g_i=g_j$ for some $i\neq j$, so
that $\ell_i$ and $\ell_j$ are within $2\delta$ of each other,
contradicting the choice of the sequence.

Finally, note that we could have taken $f_i=gf^i$ for some $g\in F$
that does not commute with $f$, and that the argument shows that
$f_i\not\sim f$ for all but finitely many $i$.
\end{proof}

\begin{thm}\label{wpd}
Suppose that $G$ and $X$ satisfy $\text{\it WPD}$.
Then $\tqh(G)$ is infinite dimensional.
\end{thm}

\begin{proof}
This is a consequence of Theorem \ref{blueprint} and Proposition
\ref{basic}.
\end{proof} 

In order to avoid passage to finite index subgroups, we will need a
slight extension of Theorem \ref{wpd}.

\begin{thm} \label{wpd2}
Suppose that $G$ and $X$ satisfy $\text{\it WPD}$. For $p\geq 1$ form the
semi-direct product $\tilde G=G^p\rtimes S_p$ where $G^p=G\times
G\times\cdots\times G$ is the $p$--fold cartesian product and $S_p$ is the
symmetric group on $p$ letters acting on $G^p$ by permuting the
factors. Let $\tilde H<\tilde G$ be any subgroup and let $H=\tilde
H\cap G^p$ (subgroup of $\tilde H$ of index $\leq p!$). Thus $H$ has
$p$ actions on $X$ obtained by projecting to various coordinates. If
at least one of these actions satisfies $\text{\it WPD}$ (equivalently, it is
nonelementary) then $\tqh(\tilde H)$ is
infinite dimensional.
\end{thm}

Note that Theorem \ref{wpd} implies that $\tqh(H)$ is
infinite dimensional.

\begin{proof} The details of this proof are similar to the proof of
Theorem \ref{wpd} and we only give a sketch. We will use the following
principle in this proof. If $F$ is a rank 2 free group and $\phi\co F\to
G$ a homomorphism then there is a rank 2 free subgroup $F'<F$ such
that either $\phi(F')$ contains no hyperbolic elements or else $\phi$
is injective on $F'$ and $\phi(F')$ is Schottky (it follows from $\text{\it WPD}$
that either $\phi(F)$ contains two independent hyperbolic elements, in
which case the latter possibility can be arranged, or $\phi(F)$
contains no hyperbolic elements, or $\phi(F)$ is virtually cyclic, and
then the first possibility holds).

Say the first projection of $H$ induces an action which is
$\text{\it WPD}$. Therefore there is a free group $F=<x,y>\subset H$ such that
the first projection of $F$ is Schottky. Now apply the above principle
with respect to each coordinate to replace $F$ by a subgroup so that
each coordinate action is either Schottky or has no hyperbolic elements. For
concreteness, we assume that coordinates $1,2,\cdots,k$ are Schottky
and $k+1,\cdots,p$ have no hyperbolic elements ($1\leq k\leq p$). We still
call $x$ and $y$ the basis elements of $F$. 

We will adopt the convention in this proof that for $f\in F$ the
$r^{th}$ projection of $f$ is denoted by $_rf$.

The proof of Proposition \ref{basic}(5) (see the last sentence) shows
that after replacing $y$ by $xy^N$ for some $N$ if necessary, we may
assume that $_rx\not\sim {_ry}$ for $r=1,2,\cdots,k$. Next, elements
$f=x^{n_1}y^{m_1}x^{k_1}y^{-l_1}$ and
$g=x^{n_2}y^{m_2}x^{k_2}y^{-l_2}$ for $0\ll n_1\ll m_1\ll k_1\ll
l_1\ll n_2\ll m_2\ll k_2\ll l_2$ will have the property that
$_rf\not\sim {_rg}^{\pm 1}$ for $r=1,2,\cdots,k$ (see Claims 1 and 2
in the proof of Proposition \ref{10}). We could then construct a
sequence $f_1,f_2,\cdots$ as in the proof of Proposition \ref{10} (in
the same manner as in the previous sentence) so that $_rf_i\not\sim
{_rf_j}^{\pm 1}$ for $i\neq j$ and $_rf_i\not\sim {_rf_i}^{-1}$. 

In
addition, we want to arrange that $_1f_j\not\sim {_rf_j}^{-1}$ for
$j=1,2,\cdots$ (note that we cannot hope to arrange $_1f_j\not\sim
{_rf_j}$ since $G$ might have the same $1^{st}$ and $r^{th}$
projections). This can be done by modifying the expression for $f_j$
so that it reads (for example)
$$f_j=x^{-s_j}y^{-t_j}x^{n_j}y^{m_j}x^{k_j}y^{-l_j}$$
with $0\ll s_1\ll t_1\ll n_1\ll m_1\ll k_1\ll l_1\ll s_2\ll\cdots$. 
The idea is that $_1f_j\sim {_rf_j}^{-1}$ would force the situation
where a long string of $_1y$'s is close to both a long string of $_rx$'s
and a long string of $_rx^{-1}$'s, implying $_rx\sim {_rx}^{-1}$. Of
course, it can be arranged that this is false by replacing $(x,y)$
with $(f_1,f_2)$ from the previous paragraph.

We now define
quasi-homomorphisms $h_i\co G^p\to\R$ by the formula
$$h_i(g_1,g_2,\cdots,g_p)=h_{({_1f_i}^{a_i})}(g_1)+
\cdots+h_{({_1f_i}^{a_i})}(g_p)$$ for large $a_i$. These maps clearly
extend to quasi-homomorphisms on $\tilde G=G^p\rtimes S_p$. The first
summand in the above formula is unbounded on the cyclic subgroup
$<f_i>$ and it is positive on large positive powers of $f_i$. The
second through $k^{th}$ summands are nonnegative on large powers of
$f_i$ thanks to the fact that $_1f_i\not\sim {_rf_i}^{-1}$ for $k\geq
r>1$. Finally, the other summands are bounded on $<f_i>$ since $_rf_i$
is not hyperbolic for $r>k$. Thus $h_i$ is unbounded on $<f_i>$. A
similar argument shows that $h_i$ is bounded on $<f_j>$ for $j<i$, so
that the elements of $QH(\tilde H)$ induced by $h_1,h_2,\cdots$ are
linearly independent. By choosing the $f_i$'s to lie in the commutator
subgroup of $G$ as before, we obtain an infinite linearly independent
set in $\tqh(F)$ and hence in $\tqh(\tilde H)$.
\end{proof}

\section{Mapping class groups}

Let $S$ be a compact orientable surface of genus $g$ and $p$
punctures. We consider the associated mapping class group $\text{\it MCG}(S)$ of
$S$. This group acts on the {\it curve complex} $X$ of $S$ defined by
Harvey \cite{harvey} and successfully used in the study of mapping
class groups by Harer \cite{harer:vcd}, \cite{harer:stability} and by
N\thinspace V Ivanov \cite{ivanov:91}, \cite{ivanov:97}. For
our purposes, we will restrict to the 1--skeleton of (the barycentric
subdivision of) Harvey's complex,
so that $X$ is a graph whose vertices are isotopy classes of
essential, nonparallel, nonperipheral, pairwise disjoint simple closed
curves in $S$ (also called {\it curve systems}) and two distinct
vertices are joined by an edge if the corresponding curve systems can
be realized simultaneously by pairwise disjoint curves. In certain
sporadic cases $X$ as defined above is 0--dimensional or empty
(this happens when there are no curve systems consisting of two
curves, ie, when $g=0$, $p\leq 4$ and when $g=1$, $p\leq 1$). In the
theorem below these cases are excluded (one could rectify the
situation by declaring that in those cases two vertices are joined by
an edge if the corresponding curves can be realized with only one
intersection point). The mapping class group
$\text{\it MCG}(S)$ acts on $X$ by $f\cdot a=f(a)$.

H Masur and Y Minsky proved the following
remarkable result.

\begin{thm}{\rm\cite{minsky-masur:cc1}}\qua
The curve complex $X$ is $\delta$--hyperbolic. An element of\break
$\text{\it MCG}(S)$
acts hyperbolically on $X$ if and only if it is \pa.\qed
\end{thm}

The following lemma is well-known (see \cite[Theorem 2.7]{cb:surfaces}).

\begin{lemma}\label{casson}
Suppose that $a$ and $b$ are two curve systems on $S$ that intersect
minimally and such that $a\cup
b$ fills $S$. Then the intersection $S(a,b)$ of the stabilizers of $a$
and of $b$
in $\text{\it MCG}(S)$ is finite.
\end{lemma}

We remark that $a\cup b$ fills $S$ if and only if $d(a,b)\geq 3$ in
the curve complex.

\begin{proof} 
Let $g$ be in the stabilizer of both $a$ and $b$. Then there is an
isotopy of $g$ so that $g(a)=a$ and $g(b)=b$. It follows that for some
$N>0$ depending only on the complexity of the graph $a\cup b$ we have
that $g^N$ is isotopic to the identity. Therefore $S(a,b)$ consists of
elements of finite order and is consequently finite (every torsion
subgroup of a finitely generated virtually torsion-free group is finite).
\end{proof}

\begin{prop}
Let $S$ be a nonsporadic surface.
The action of $\text{\it MCG}(S)$ on the curve complex $X$ satisfies $\text{\it WPD}$.
\end{prop}

\begin{proof}
The first two bullets in the definition of $\text{\it WPD}$ are clear. Our proof
of the remaining property is motivated by Feng Luo's proof (as explained
in \cite{minsky-masur:cc1}) that the
curve complex has infinite diameter. We recall
the construction and the basic properties of Thurston's space of
projective measured foliations on $S$ (see \cite{thurston:surfaces} and
\cite{flp:surfaces}). Let $\C$ be the set of all curve systems
in $S$ and by
$$I\co \C\times\C\to [0,\infty)$$ denote the intersection pairing,
ie, $I(a,b)$ is the smallest number of intersection points between
$a$ and $b$ after a possible isotopy. Let $\r=(0,\infty)$ and by
$\r\C$ denote the space of formal products $ta$ for $t\in\r$ and $a\in
\C$ where we identify $\C$ with the subset $1\C$. Extend $I$ to
$\r\C\times \r\C$ by $$I(ta,sb)=tsI(a,b).$$ Consider the associated
function $$J\co \r\C\to [0,\infty)^\C$$ defined by $$J(ta)=(sb\mapsto
I(ta,sb)).$$ Then $J$ is injective and we let $\mf$ denote the closure
of the image of $J$. An element of $\mf$ can be viewed as a measured
foliation on $S$. The pairing $I$ extends to a continuous function
$$I\co \mf\times\mf\to [0,\infty).$$

There is a natural action of $\r$ on $\mf$ given by scaling. The orbit
space $\pmf$ is Thurston's space of projective measured foliations
and it is homeomorphic to the sphere of dimension $6g+2p-7$ (assuming
this number is positive). The intersection pairing is not defined on
$\pmf\times\pmf$ but note that the statement $I(\Lambda,\Lambda')=0$
makes sense for $\Lambda,\Lambda'\in\pmf$. The mapping class group
$\text{\it MCG}(S)$ of $S$ acts on $\C$ by $f\cdot a=f(a)$ and there is an
induced action on $\r\C$, $\mf$, and $\pmf$.

Let $f\in \text{\it MCG}(S)$ be a \pa mapping class. Then $f$ fixes exactly two
points in $\pmf$ and one point $\Lambda_+$ is attracting while the
other $\Lambda_-$ is repelling. All other points converge to
$\Lambda_+$ under forward iteration and to $\Lambda_-$ under backward
iteration. It is known that $I(\Lambda_+,\Lambda)=0$ implies
$\Lambda=\Lambda_+$ and similarly for $\Lambda_-$. 
Continuity of $I$ implies the following fact: 

{\it If $U$ is a neighborhood
of $\Lambda_+$ then there is a neighborhood $V$ of $\Lambda_+$ such
that if $\Lambda,\Lambda'\in\pmf$, $I(\Lambda,\Lambda')=0$ and
$\Lambda\in V$ then $\Lambda'\in U$.}

We will use the terminology that $V$ is {\it adequate} for $U$ if the
above sentence holds. A similar fact (and terminology) holds for
neighborhoods of $\Lambda_-$.

Given $C>0$, choose closed neighborhoods $U_0\supset U_1\supset
U_2\supset\cdots\supset U_N$ of $\Lambda_+$ and  $V_0\supset V_1\supset
V_2\supset\cdots\supset V_N$ of $\Lambda_-$ with $N>C$ so that
\begin{itemize}
\item $U_{i+1}$ is adequate for $U_i$ and $V_{i+1}$ is adequate for
$V_i$, and
\item if $\Lambda\in U_0$ and $\Lambda'\in V_0$ then
$I(\Lambda,\Lambda')\neq 0$.
\end{itemize}
Assume now that two curve systems $a$ and $b$ belong to a quasi-axis $\ell$
of a \pa mapping class $f$ and that they are sufficiently far away from
each other, so that after applying a power of $f$ and possibly
interchanging $a$ and $b$ we may assume that
$a\in U_N$ and $b\in V_N$. Assume, by way of contradiction, that $g_n$
is an infinite sequence in $\text{\it MCG}(S)$ and $d(a,g_n(a))\leq N$,
$d(b,g_n(b))\leq N$ for all $n$. Note that if $c$ is a curve system
disjoint from $a$ then $c\in U_{N-1}$, and inductively if $d(a,c)\leq
N$ then $c\in U_0$. We therefore conclude that $g_n(a)\in U_0$ and
$g_n(b)\in V_0$. After passing to a subsequence, we may assume that
the sequence $g_n(a)$ converges to $A\in U_0$ and $g_n(b)\to B\in
V_0$. Note that $I(A,B)\neq 0$ by the choice of $U_0$ and $V_0$.

First suppose that the curve systems $g_n(a)$ are all different. To
obtain convergence in $\mf$ one is required to first rescale by some
$r_n>0$, ie, $\frac1{r_n}g_n(a)\to \tilde A\in\mf$ where $r_n$ can be
taken to be the length of $g_n(a)$ in some fixed hyperbolic
structure on $S$. Under the assumption that $g_n(a)$ are all distinct, we see
that $r_n\to\infty$ and this implies that $I(\tilde A,\tilde B)=0$
ie, that $I(A,B)=0$, contradiction.
The case when $g_n(b)$ are all distinct is similar.

Finally, if $g_n(a)$ and $g_n(b)$ take only finitely many values, we
may assume by passing to a subsequence that both $g_n(a)$ and $g_n(b)$
are constant. But then $g_n^{-1}g_m\in S(a,b)$ and Lemma \ref{casson}
implies that the sequence $g_n$ is finite.
\end{proof}

The following is the main theorem in this note. H Endo and
D Kotschick \cite{e-k:bc} have shown using 4--manifold topology and
Seiberg--Witten invariants that $\tqh(\text{\it MCG}(S))\neq 0$ when $S$ is
hyperbolic. M Korkmaz \cite{korkmaz:bc} also proved\break $\tqh(\text{\it MCG}(S))\neq
0$, and in addition that $\tqh(\text{\it MCG}(S))$ is infinite dimensional when
$S$ has low genus. The nontriviality, and even
infinite-dimensionality, of $\tqh(\text{\it MCG}(S))$ was conjectured by Morita
\cite[Conjecture 6.19]{morita:survey}.

\begin{thm}\label{main}
Let $G$ be a subgroup
of $\text{\it MCG}(S)$ which is not virtually abelian.
Then
$\dim\tqh(G)=\infty$.
\end{thm}

\proof
We first deal with the sporadic cases. When $g=0$, $p\leq 3$ and when
$g=1$, $p=0$ the mapping class group $\text{\it MCG}(S)$ is finite. When $g=0$,
$p=4$ and when $g=p=1$ then $\text{\it MCG}(S)$ is word hyperbolic (in fact, the
quotient by the finite center is virtually free) and instead of
considering the action on a curve complex we can look at the action on
the Cayley graph. This action is properly discontinuous and therefore
the restriction to any subgroup which is not virtually cyclic
satisfies $\text{\it WPD}$. The statement then follows from Theorem \ref{wpd}.

Now we assume that $S$ is not sporadic.  By the classification of
subgroups (see \cite[Theorem 4.6]{mp:subgroups}) there are 4 cases.
\begin{itemize}
\item $G$ contains two independent \pa homeomorphisms. Then the action
of $G$ on the curve complex $X$ for $S$ satisfies the assumptions of
Theorem \ref{wpd} so $\tqh(G)$ is infinite dimensional.
\item $G$ fixes a pair $\Lambda_{\pm}$ of foliations corresponding to
a \pa homeomorphism. Then $G$ is virtually cyclic.
\item $G$ is finite.
\item There is a curve system $c$ on $S$ invariant under $G$. Choose
$c$ to be maximal possible and cut
$S$ open along $c$. Consider the mapping class group of the cut
open surface $S'$ where we collapse each boundary component to a
puncture. Since $c$ is maximal and $G$ is not virtually abelian, there
is an orbit $S_1',S_2',\cdots,S_p'$ of components of $S'$ and there is
a subgroup of $G$ that preserves these components and whose
restriction to a component contains two independent \pa
homeomorphisms. Pass to the quotient of $G$ corresponding to the
restriction to $S_1'\cup\cdots\cup S_p'$. The mapping class group of
$S_1'\cup\cdots\cup S_p'$ can be identified with $\text{\it MCG}(S_1')^p\rtimes
S_p$ so the result in this case follows from Theorem \ref{wpd2}.\qed
\end{itemize}

\medskip

The following is a version of superrigidity for mapping class
groups. It was conjectured by N\thinspace V Ivanov and proved by Kaimanovich
and Masur \cite{km:poisson} in the case when the image group contains
independent \pa homeomorphisms and it was extended to the general case
by Farb and Masur \cite{fm:rigidity} using the classification of
subgroups of $\text{\it MCG}(S)$ as above. Our proof is different in that it does
not use random walks on mapping class groups, but instead uses the
work of M Burger and N Monod \cite{bm:lattices} on bounded
cohomology of lattices. Note also that for this application we only
need a weak version of our result, namely that $\tqh(G)\neq 0$ when
$G\subset \text{\it MCG}(S)$ is not virtually abelian.

\begin{cor}
Let $\Gamma$ be an irreducible lattice in a connected semi-simple Lie
group $G$ with no compact factors, with finite center, and of rank
$>1$. Then every homomorphism $\Gamma\to \text{\it MCG}(S)$ has finite image.
\end{cor}

\begin{proof} Let 
$\phi\co \Gamma\to \text{\it MCG}(S)$ be a homomorphism. 
By the Margulis--Kazhdan theorem \cite[Theorem 8.1.2]{zimmer:book}
either the image of $\phi$ is finite or the kernel of $\phi$ is
contained in the center. When $\Gamma$ is a nonuniform lattice, the
proof is easier and was known to Ivanov before the work of
Kaimanovich--Masur (see Ivanov's comments to Problem 2.15 on Kirby's
list). Since the rank is $\geq 2$ the lattice $\Gamma$ then contains a
solvable subgroup $N$ which does not become abelian after quotienting
out a finite normal subgroup. If the kernel is finite, then $\phi(N)$
is a solvable subgroup of $\text{\it MCG}(S)$ which is not virtually abelian,
contradicting \cite{blm:s2a}.

Now assume that $\Gamma$ is a uniform lattice. If the kernel $Ker(\phi)$ is
finite then there is an unbounded quasi-homomorphism
$q\co Im(\phi)\to\R$ by Theorem \ref{main}. But then $q\phi\co \Gamma\to\R$
is an unbounded quasi-homomorphism contradicting the Burger--Monod
result that every quasi-homomorphism $\Gamma\to\R$ is bounded. 
\end{proof}

\begin{remark*}
When the center $Z(G)$ of $G$ is infinite, one can show that every
homomorphism $\phi\co \Gamma\to \text{\it MCG}(S)$ has virtually abelian image, as
follows. The key is that in this case the intersection $\Gamma\cap
Z(G)$ has finite index in $Z(G)$ and the projection of $\Gamma$ in
$G/Z(G)$ is a lattice $L$ in $G/Z(G)$, which is a Lie group of rank
$>1$. Choose $g\in\Gamma\cap Z(G)$ of infinite order such that the
Nielsen--Thurston representative of $\phi(g)$ has no rotations
(ie, the canonical invariant curve system cuts the surface into
invariant subsurfaces and on each subsurface (with boundary components
collapsed to punctures) $\phi(g)$ is isotopic to a pseudo-Anosov
homeomorphism or to the identity). There is an induced map
$\phi\co \Gamma\to C(\phi(g))$, the centralizer of $\phi(g)$. Moreover,
from the Nielsen--Thurston theory we have a homomorphism
$\psi\co C(\phi(g))\to \text{\it MCG}(S_1)\times \text{\it MCG}(S_2)\times\cdots\times
\text{\it MCG}(S_k)$ given by restricting to the subsurfaces on which $\phi(g)$
is identity. The kernel of $\psi$ is virtually abelian. So it suffices
to argue that the image of the composition $\psi\phi\co \Gamma/<g>\to
\text{\it MCG}(S_1)\times \text{\it MCG}(S_2)\times\cdots\times \text{\it MCG}(S_k)$ is virtually
abelian. Now $\Gamma/<g>$ is a lattice in $G/<g>$, a Lie group of rank
$>1$ and the rank of the center $Z(G/<g>)=Z(G)/<g>$ is smaller than
the rank of $Z(G)$. So the statement follows by induction on the rank
of the center.
\end{remark*}

\begin{remark*}
It can also be shown that the image of $\Gamma\to \text{\it MCG}(S)$ must be
finite, even if $Z(G)$ has infinite center. If the image is infinite,
we can assume (by passing to a subgroup of $\Gamma$ of finite index)
that it is torsion-free and abelian of finite rank. If $G$ and
$\Gamma$ satisfy Kazhdan's property (T), then the abelianization of
$\Gamma$ is finite and this is the case when $G$ has no rank 1
factors. It is also true that the abelianization of $\Gamma$ is finite
when $G$ is any higher rank group and this follows from a deep work of
Prasad--Raghunathan and Deligne.
\end{remark*}

\begin{remark*}
Theorem \ref{main} also implies that S--arithmetic groups and certain
groups of automorphisms of trees do not occur as subgroups of mapping
class groups. These are the groups $\Gamma$ for which Burger--Monod
show $\tqh(\Gamma)=0$ \cite{bm:lattices}, \cite{burger-monod2},
\cite{burger-monod3}.
\end{remark*}

{\bf Acknowledgments}\qua The authors would like to thank Yair Minsky and Gopal
Prasad for useful discussions. The first author
  gratefully acknowledges the support by the National Science
  Foundation. The second author appreciates the hospitality of the
  mathematics department of the University of Utah.

\end{document}